\documentclass{amsart}
\usepackage{amsaddr}
\usepackage{amssymb,amsmath,amsfonts}
\usepackage{latexsym,amscd}
\usepackage{pgf,tikz,pgfplots}
\pgfplotsset{compat=1.14}
\usepackage{mathrsfs}
\usetikzlibrary{arrows,backgrounds,calc,decorations,patterns}
\usepackage{url}
\usepackage[numbers]{natbib}

\usetikzlibrary {patterns.meta}
\pgfdeclarepattern{
  name=hatch,
  parameters={\hatchsize,\hatchangle,\hatchlinewidth},
  bottom left={\pgfpoint{-.1pt}{-.1pt}},
  top right={\pgfpoint{\hatchsize+.1pt}{\hatchsize+.1pt}},
  tile size={\pgfpoint{\hatchsize}{\hatchsize}},
  tile transformation={\pgftransformrotate{\hatchangle}},
  code={
    \pgfsetlinewidth{\hatchlinewidth}
    \pgfpathmoveto{\pgfpoint{-.1pt}{-.1pt}}
    \pgfpathlineto{\pgfpoint{\hatchsize+.1pt}{\hatchsize+.1pt}}
    \pgfpathmoveto{\pgfpoint{-.1pt}{\hatchsize+.1pt}}
    \pgfpathlineto{\pgfpoint{\hatchsize+.1pt}{-.1pt}}
    \pgfusepath{stroke}
  }
}

\tikzset{
  hatch size/.store in=\hatchsize,
  hatch angle/.store in=\hatchangle,
  hatch line width/.store in=\hatchlinewidth,
  hatch size=5pt,
  hatch angle=0pt,
  hatch line width=.5pt,
}

\newtheorem{thm}{Theorem}
\newtheorem{lem}[thm]{Lemma}

\newtheorem{prop}[thm]{Proposition}
\newtheorem{cor}[thm]{Corollary}

\theoremstyle{definition}

\newtheorem{defn}[thm]{Definition}
\newtheorem*{fact}{Fact}

\newtheorem{example}{Example}

\def\R{\mbox{{\bf R}}}

\def\Z{\mbox{{$\mathbb{Z}$}}}

\def\cp{G}

\DeclareMathOperator{\Area}{Area}

\DeclareMathOperator{\word}{Word}

\newcommand{\s}{\sigma}

\def \sb #1 {\overline{\s}_{#1}}

\def\aa #1 {\textcolor{blue}{#1}}

\begin{document}
\title[Non-abelian Sperner]{Integer Area Dissections of Lattice Polygons\\
via a Non-Abelian Sperner's Lemma}
%
\author{Aaron Abrams}
\address{School of Data Science and Mathematics Department, University of Virginia\\ Mathematics Department, Washington and Lee University}
\author{Jamie Pommersheim}
\address{Mathematics and Statistics Department, Reed College}

\email{abrams.aaron@gmail.com, jamie@reed.edu}


\begin{abstract}
    We give a simple and complete description of those convex lattice polygons in the plane that can be dissected into lattice triangles of integer area.  A new version of Sperner's Lemma plays a central role.
\end{abstract}
\keywords{Dissection, lattice polygon, Sperner's Lemma}

\subjclass{52B45,52C05}

\maketitle

\section{Prologue}

    This article addresses the following main question:  
    \begin{quote}
        \emph{Which convex lattice polygons can be dissected into lattice triangles of area 1?}
    \end{quote}
    
    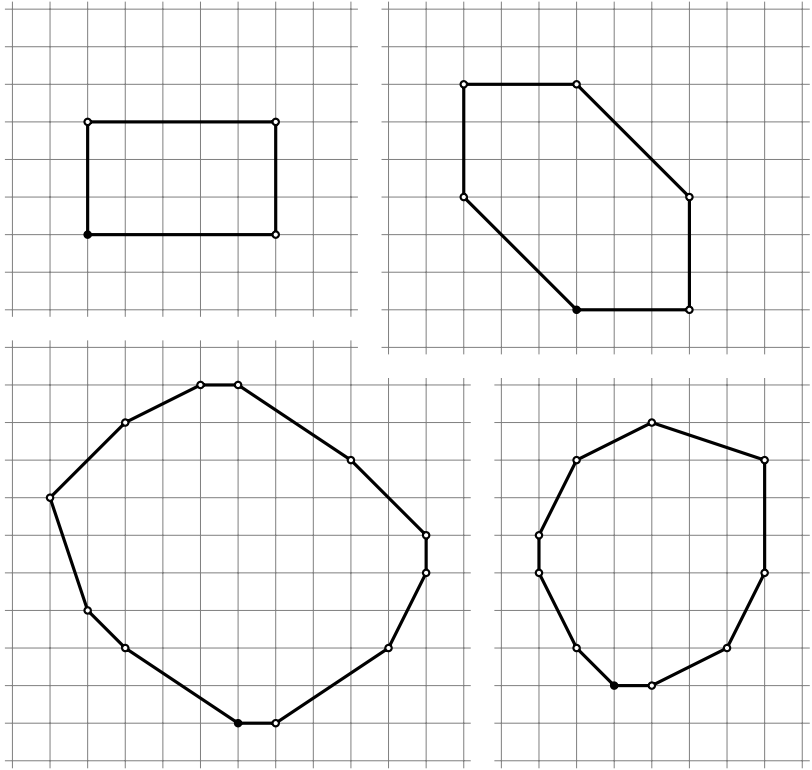
\begin{figure}
        \centering
        \begin{tikzpicture}[scale=.5]
        \input ./graphpaper.tex
        \end{tikzpicture}
        \caption{Which of these polygons can be cut into lattice triangles of area 1?}
        \label{fig:polygons}
    \end{figure}

    A \emph{lattice polygon} is a polygon $P$ in the plane whose vertices have integer coordinates.
    See Figure \ref{fig:polygons}.  To \emph{dissect} a polygon into triangles is to express the polygon as a union of triangles whose interiors do not overlap, as exemplified in Figure \ref{fig:dissection}.
    
    The area of a lattice polygon is always half an integer, and in fact every lattice polygon can be dissected into lattice triangles of area $1/2$. (This important fact features prominently in many of the proofs of Pick's Theorem for the area of a lattice polygon \cite{thebook, nz, varberg, gs}.)
    
    A lattice polygon with integer area, however, cannot necessarily be dissected into triangles of area 1. This is demonstrated by the smallest possible example, a $1\times 1$ square.
    
    In this article we will establish an easily checkable necessary and sufficient condition for a convex lattice polygon to have a dissection into lattice triangles of area 1.
    
    Although this problem presents itself as a mix of combinatorial geometry and number theory,
    it ultimately yields to an argument that is inherently topological.
    The hero of the story is a new version of Sperner's Lemma, 
    which we will present in Section \ref{sec:nonabelian}.

\section{Setup and Main Theorem}\label{sec:intro}

    As we try to get a feel for how to think about the main question, a few outer layers can be peeled off fairly readily.

    \subsection*{Treating Triangles}
        Suppose the whole polygon $P$ is a triangle of integer area, or that we have dissected $P$ into lattice triangles of integer area.
        In this case, it is \emph{always} possible to finish the job:  
        any such triangle can be dissected into lattice triangles of area $1$. 
        We prove this in Section \ref{sec:appendix}, using a bit of elementary number theory and linear algebra that is independent of the rest of the paper.
        
        As a consequence, we can reformulate our main question as follows.
        \begin{quote}
            \emph{Which convex lattice polygons can be dissected into lattice triangles of integer area?}
        \end{quote}
        This also motivates the following definition.
        \begin{defn}[Integral dissection]  
            A dissection of a lattice polygon into lattice triangles of integer area is called an \emph{integral dissection}. 
        \end{defn}

    \begin{figure}
    \centering
    \begin{tikzpicture}[scale=.5]
    \draw[step=1, gray, very thin] (-.2,-.2) grid (7.2,6.2);

\coordinate (P) at (2,1);
\draw[thick] (P) -- ++(3,0)
    -- ++(1,3)
    -- ++(-3,1)
    -- ++(-2,-2)
    -- cycle;

\draw[thick] (P)++(-1,2) -- +(4,0); 
\draw[thick] (P)++(0,2) -- ++(0,-2) -- ++(1,2) -- ++(0,-2); 
\draw[thick] (P)++(1,2) -- ++(2,-2) -- ++(0,2) -- ++(-2,2); 
\draw[thick] (P)++(3,2) -- ++(1,1);

\draw[fill=white] (P) circle (2pt);
\draw[fill=white] (P)+(3,0) circle (2pt);
\draw[fill=white] (P)+(4,3) circle (2pt);
\draw[fill=white] (P)+(1,4) circle (2pt);
\draw[fill=white] (P)+(-1,2) circle (2pt);

\draw[fill=white] (P)+(0,2) circle (2pt);
\draw[fill=white] (P)+(1,2) circle (2pt);
\draw[fill=white] (P)+(3,2) circle (2pt);
\draw[fill=white] (P)+(1,0) circle (2pt);
    \end{tikzpicture}
    \caption{A dissection of a lattice pentagon into integer area lattice triangles.}
    \label{fig:dissection}
    \end{figure}
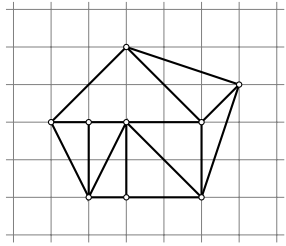

    \subsection*{Pursuing Parity}
        Parity considerations play a prominent role in our story, as one might expect. The \emph{parity} of a lattice point is one of the four pairs 
        \[ A=(\text{even}, \text{even}),\ B=(\text{odd},\text{even}),\ C=(\text{odd}, \text{odd}),\ D=(\text{even}, \text{odd}), \]
        according to the parities of its coordinates.
        We will often refer to the parities $A,B,C,D$ as \emph{colors} and speak of lattice points as being \emph{colored by parity}. 

        Notice that two lattice points have the same color if and only if the midpoint of the segment joining them is a lattice point. Thus if $P$ is a lattice triangle with two vertices of the same color, we can split $P$ into two lattice triangles of equal area, so $P$ itself has integer area.  
        Conveniently, the converse also holds:  a lattice triangle has integer area if and only if two of its vertices have the same color.  We prove this Parity Proposition in Section \ref{sec:triangles}.  (A similar fact appears as Lemma 1 in both \cite{thomas} and \cite{jz}.)
        This leads to yet a third formulation of the main question.
        \begin{quote}
            \emph{Which convex lattice polygons can be dissected into lattice triangles such that each triangle has (at least) two vertices of the same color?}
        \end{quote}
        Readers familiar with Sperner's Lemma may already sense its presence here.

    \subsection*{Dissecting Diagonally}  
        One way to dissect a convex polygon $P$ is to draw diagonals, without adding new vertices, until $P$ is dissected into triangles.  This special type of dissection in which no new vertices are introduced is called a \emph{diagonal dissection}.  
                
        Not every diagonal dissection is an integral dissection.  But if we seek an integral diagonal dissection, then by the Parity Proposition it suffices to know the colors of the vertices of $P$.

        If $P$ is a lattice polygon, then the colors of the vertices of $P$, read counterclockwise, form a \emph{cyclic word}, by which we mean a finite string defined up to cyclic permutation.  We use notation $(X_1X_2\ldots X_n)$ to denote a cyclic word; this is the same cyclic word as $(X_2X_3\ldots X_nX_1)$, and so on.
        We treat the subscripts as integers modulo $n$, so for instance the notation $X_{n+1}$ refers to the symbol $X_1$.
                        
        \begin{defn}[Boundary Word]
            Let $P$ be a lattice polygon, with vertices colored according to parity.
            The {\it boundary word} of $P$ is the cyclic word $\word(P)$ obtained by recording the colors of the vertices in counterclockwise order around $P$.
        \end{defn}
        
        For example, each of the polygons in Figure \ref{fig:polygons} has one vertex distinguished by a solid dot, which represents the origin.  One sees that the rectangle in the upper left has boundary word $(ABCD)$. The decagon $P$ in the lower right has $\word(P)=(ABCDACBADC)$.

        \begin{example}\label{ex:diagdiss} 
            Consider the dodecagon $P$ in the lower left of Figure \ref{fig:polygons}.  This dodecagon has boundary word $(ABABCCDCBBDB)$. An integral diagonal dissection of $P$ is illustrated in Figure \ref{fig:diagonal}. We have drawn $P$ here as if it were regular in order to emphasize that what matters for this process is only the coloring of the vertices, and not their exact geometric positions.
            
            To construct this dissection, we might start by drawing a diagonal joining the two $A$'s, creating a triangle labeled $ABA$ that has integer area by the Parity Proposition. 
            Eventually we might end up with the picture on the right, which shows (the combinatorial structure of) an integral diagonal dissection of $P$.  
            
            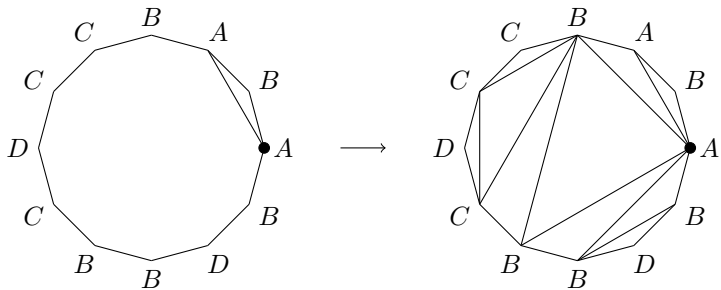
\begin{figure}[h]
                \centering
                \begin{tikzpicture}[scale=.5]
                \foreach \x in {0,...,12}
    \coordinate (\x) at (30*\x : 3);
    


\matrix[row sep = .5cm, column sep = .5cm]{
    \foreach \x/\y in {0/1,1/2,2/3,3/4,4/5,5/6,6/7,7/8,8/9,9/10,10/11,11/12}
    \draw (\x) -- (\y);
    \draw[fill=black] (0) circle (2pt);
    \draw (0) node [anchor=180]{$A$};
    \draw (1) node [anchor=210]{$B$};
    \draw (2) node [anchor=240]{$A$};
    \draw (3) node [anchor=270]{$B$};
    \draw (4) node [anchor=300]{$C$};
    \draw (5) node [anchor=330]{$C$};
    \draw (6) node [anchor=0]{$D$};
    \draw (7) node [anchor=30]{$C$};
    \draw (8) node [anchor=60]{$B$};
    \draw (9) node [anchor=90]{$B$};
    \draw (10) node [anchor=120]{$D$};
    \draw (11) node [anchor=150]{$B$};
    \draw (0)--(2);
    &
    \draw[->] (0,0)--(.6,0);
    &
    \foreach \x/\y in {0/1,1/2,2/3,3/4,4/5,5/6,6/7,7/8,8/9,9/10,10/11,11/12}
    \draw (\x) -- (\y);
    \draw[fill=black] (0) circle (2pt);
    \draw (0) node [anchor=180]{$A$};
    \draw (1) node [anchor=210]{$B$};
    \draw (2) node [anchor=240]{$A$};
    \draw (3) node [anchor=270]{$B$};
    \draw (4) node [anchor=300]{$C$};
    \draw (5) node [anchor=330]{$C$};
    \draw (6) node [anchor=0]{$D$};
    \draw (7) node [anchor=30]{$C$};
    \draw (8) node [anchor=60]{$B$};
    \draw (9) node [anchor=90]{$B$};
    \draw (10) node [anchor=120]{$D$};
    \draw (11) node [anchor=150]{$B$};
    \draw (0)--(2);
    \draw (0)--(3);
    \draw (3)--(5);
    \draw (5)--(7);
    \draw (3)--(7);
    \draw (3)--(8);
    \draw (0)--(8);
    \draw (0)--(9);
    \draw (9)--(11);
    \\
    %
    %
    };
                \end{tikzpicture}
                \caption{A diagonally dissected dodecagon.}
                \label{fig:diagonal}
            \end{figure}
        \end{example}

    \subsection*{Conjuring Contraction} 
        We will see that integral diagonal dissections can be interpreted, and discovered, using only the boundary word.
        Consider the following operation on cyclic words.
        
        \begin{defn}[contracting step, contractible] \label{def:contract}
            Let $w$ be a cyclic word $(X_1 X_2 \dots X_n)$.       
            As usual indices are taken modulo $n$, and we assume $n\ge 2$.
            \begin{quote}
                \emph{Contracting step}:\   If letters $X_{i-1},X_i,X_{i+1}$ are not all distinct, then delete $X_i$ from $w$.
            \end{quote}
                The cyclic word $w$ is \emph{contractible} if there is a sequence of contracting steps starting with $w$ that results in a word $(X)$ consisting of a single letter. 
        \end{defn}

        Note that a contracting step can be used to reduce either of the words $(XY)$ or $(XX)$ to the word $(X)$, so a word is contractible if and only if it can be reduced until its length is less than three.

        Contractibility of $\word(P)$ is equivalent to the existence of an integral diagonal dissection of $P$,
        as we prove in Section \ref{sec:3iff4}.  In Section \ref{sec:alg} we show that it is algorithmically easy to decide whether a given cyclic word is contractible.  Together, these results solve the problem of finding integral diagonal dissections of convex lattice polygons, which is great progress towards our main question.

    \subsection*{Main Theorem}
        We have just asserted that if the boundary word of $P$ is contractible, then there is a dissection (in fact, a diagonal dissection) of $P$ into lattice triangles of integer area.  As it turns out, the converse is also true!
        That is, if the boundary word of $P$ is not contractible, then not only is there no integral \emph{diagonal} dissection of $P$ (again, as asserted above), but there is no integral dissection of $P$ at all.
        This is the heart of our main theorem.
        \begin{thm}\label{thm:stacked}
            For any convex lattice polygon $P$, the following are equivalent.
            \begin{enumerate} 
            \item[(1)]
            $P$ can be dissected into lattice triangles of area 1;
            \item[(2)]
            $P$ can be dissected into lattice triangles of integer area;
            \item[(3)]
            $P$ can be diagonally dissected into lattice triangles of integer area;
            \item[(4)] The boundary word of $P$ is contractible.
            \end{enumerate}
                
            Moreover, there is a linear-time algorithm to decide whether or not (1)-(4) hold.
        \end{thm}

        The final statement of Theorem \ref{thm:stacked} is proved as Proposition \ref{prop:speedy}.  Theorem \ref{thm:trianglecase} establishes the equivalence of (1) and (2), and Corollary \ref{cor:readmylips} establishes the equivalence of (3) and (4).  That (3) implies (2) is trivial, but its converse is not.  This article's \emph{raison d'\^etre} is to highlight a new ``non-abelian'' version of Sperner's Lemma (Theorem \ref{thm:nonabelian}), which is the key to showing the equivalence of (2) and (4) (and hence (3)).  
        We recommend as an exercise to try to prove that (2) implies (3) for quadrilaterals.

    \subsection*{Topology}  
        
        The resolution of our main question relies crucially on the famous Sperner's Lemma \cite{sperner}, which has been written about extensively.  
        The present paper builds on ideas of John Thomas \cite{thomas}, 
        whose work illuminated, perhaps for the first time, a beautiful and surprising connection between dissection problems and Sperner’s Lemma.  What makes the connection surprising is that questions that are highly sensitive to geometric details are resolved by an argument that is inherently topological.  

        Thomas essentially gave a complete solution to our main question for rectangles, although his motivation and context were different.  
        His reasoning can be extended to more polygons than just rectangles, but it has limitations; we discuss this further in Section \ref{sec:impossibility}.

        We were therefore delighted to find that the rest of the solution also comes from Sperner's Lemma, 
        after adapting it in a very natural way.
        Finding the right modification of Sperner's Lemma requires us to think of the lemma from a topological perspective. 
        We will end up with what we call a ``non-abelian'' version of Sperner's Lemma which will help us completely solve our main question.
        
        In the end, our geometric theorems have entirely combinatorial proofs, with no topological prerequisites.  
        Yet the topological perspective provides the theorems with extra richness, inspiration, and simplicity.

\section{Triangles} \label{sec:triangles}
    
    Let $P$ be a triangle in the plane, with vertices $(x_i, y_i)$, $i=1,2,3$. 
    We express the area of this triangle in terms of a determinant:
    \begin{equation}\label{eq:area}
    \Area(P)=\frac 12
            \left|
            \begin{matrix}
            1 & 1 & 1 \\
            x_1 & x_2 & x_3 \\
            y_1 & y_2 & y_3 
            \end{matrix}
            \right|.
    \end{equation}
    This is really a \emph{signed area}: if the points are ordered in a counterclockwise manner, the value will be positive and equal to the usual area.  Points ordered clockwise will yield a negative area.  The area is zero if and only if the three points are collinear.
    
    If $P$ is a lattice triangle, then this area is evidently one half times an integer. For the area to be an integer, the above determinant needs to be even, so the parities of the coordinates come into play.

    \begin{quote}
        \normalsize\emph{From now on, we will always assume lattice points are colored according to parity, using colors $A,B,C,D$ as defined in Section \ref{sec:intro}.  In particular, all vertices of lattice polygons and lattice dissections are automatically colored.\\}
    \end{quote}
    
    \normalsize
    \begin{prop}[Parity Proposition] \label{prop:parity}
        The area of a lattice triangle $P$ is an integer if and only if two if its vertices have the same color.
    \end{prop}
    
    \begin{proof}
        The area of $P$ is an integer if and only if the determinant in Equation (\ref{eq:area}) is an even integer.  This happens if and only if the columns of the matrix are linearly dependent modulo 2. Because of the top row, none of the columns is the zero vector mod 2 and no two columns add up to the remaining column. Thus the only possible linear dependence is for two columns to be equal mod 2. Hence the determinant is 0 mod 2 if and only if two vertices have the same color. 
    \end{proof}

    We mention a useful byproduct of this, which one could also prove directly, as it will come in handy later.
    
    \begin{cor}[Collinear Corollary]\label{cor:collinear}
        Three collinear lattice points in the plane cannot all be different colors.
    \end{cor}
    
    This is a special case of the lemma because three collinear lattice points form a triangle of area zero.

\section{Combinatorial polygons}

    As suggested by Example \ref{ex:diagdiss}, we will benefit from having a purely combinatorial description of polygons and their triangulations.  
    
    For $n\ge 3$, a \emph{combinatorial $n$-gon} is an abstract cycle graph $\cp$ with $n$ vertices; a \emph{combinatorial polygon} is a combinatorial $n$-gon for some $n$.   In contrast with lattice polygons, a combinatorial polygon does not live in the plane; it is a combinatorial object, not a geometric one.

    What does it mean to triangulate a combinatorial polygon?  
    The idea is to view a triangle as being defined by its vertices. Any 3-element set is a triangle, and an \emph{abstract triangulation} is an arbitrary collection of triangles.  
        
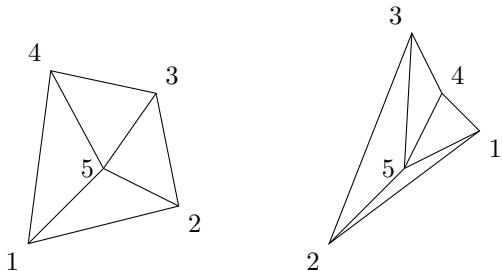
\begin{figure}[h]
        \centering
        \begin{tikzpicture}[scale=1]
        \coordinate (P) at (0,0);
\coordinate (Q) at (2,.5);
\coordinate (R) at (1.7,2);
\coordinate (S) at (.3,2.3);
\coordinate (T) at (1,1);

\draw (P)--(Q)--(R)--(S)--cycle;
\draw (P)--(T);
\draw (Q)--(T);
\draw (R)--(T);
\draw (S)--(T);

\draw (P) node[anchor=north east]{$1$};
\draw (Q) node[anchor=north west]{$2$};
\draw (R) node[anchor=south west]{$3$};
\draw (S) node[anchor=south east]{$4$};
\draw (T) node[anchor=east]{$5$};

\coordinate (P') at (4,0);
\coordinate (Q') at (6,1.5);
\coordinate (R') at (5.5,2);
\coordinate (S') at (5.1,2.8);
\coordinate (T') at (5,1);

\draw (P')--(Q')--(R')--(S')--cycle;
\draw (P')--(T');
\draw (Q')--(T');
\draw (R')--(T');
\draw (S')--(T');

\draw (P') node[anchor=north east]{$2$};
\draw (Q') node[anchor=north west]{$1$};
\draw (R') node[anchor=south west]{$4$};
\draw (S') node[anchor=south east]{$3$};
\draw (T') node[anchor=east]{$5$};
        \end{tikzpicture}
        \caption{Possible pictures of the triangulation $T$.} 
        \label{fig:abstract}
    \end{figure}

    If you wanted to convey to a friend over the phone the combinatorial structure of a triangulation, but not its specific geometric shape, this would be a good way to do it.  For example one possible triangulation $T$ has triangles $\{1,2,5\}$, $\{2,3,5\}$, $\{3,4,5\}$, $\{4,1,5\}$.  Your friend will undoubtedly picture this triangulation by drawing the vertices in the plane and drawing the triangles as in Figure \ref{fig:abstract}.  (You should agree ahead of time to draw each vertex and each edge\footnote{An edge is a 2-element subset of a triangle.} only once.)  Drawings may differ slightly, but any such picture reveals the topological structure of the triangulation:  there are four triangles that fit together around a central vertex (number 5), forming a quadrilateral with the other vertices as corners.  See Figure \ref{fig:abstract}. 
    
    Some readers will recognize these as examples of simplicial complexes.
    To focus on the topological features of the triangulation $T$, we can imagine building an abstract topological space by gluing together triangles in the pattern specified by $T$.  If two triangles share two vertices, they should be glued along the corresponding edge.  The result is called the \emph{underlying space} of $T$ and is denoted $|T|$.  It is not a subset of the plane; it is its own space.
    In general, it may or may not be possible to realize $|T|$ as a subset of the plane.  For example, the triangulation with triangles $\{1,2,3\},\{1,2,4\},\{1,2,5\}$ is not planar.
            
    Since we are ultimately interested in triangulating planar polygons, we will require our triangles to fit together to form a disk.  In this case, $|T|$ has a boundary circle which naturally inherits the structure of a combinatorial polygon.
    Figure \ref{fig:abstract}, for example, depicts a triangulation of a combinatorial $4$-gon.
    
    \begin{defn}[Triangulation] \label{def:triangulation}
        A \emph{triangulation} of a combinatorial polygon $\cp$ consists of a collection $T$ of 3-element sets called \emph{triangles}, such that the underlying space $|T|$ built out of triangles is a topological disk with boundary equal to $\cp$.  The boundary vertices are called \emph{corners}.
    \end{defn}

    Because the triangles always fit together to form a disk, every triangulation can be drawn in the plane, and in practice we infer the combinatorial structure from such a drawing.  

    The abstraction of triangulations of combinatorial polygons begins to pay off in Corollary \ref{cor:readmylips} in the next section.
    After developing this perspective a bit further, we will explicate the full relevance of these objects to the study of dissections of lattice polygons in Section \ref{sec:poofing}.

\section{Integral diagonal dissections} \label{sec:3iff4}
    Part of the main theorem asserts that our notion of contractibility for $\word(P)$ corresponds exactly to the existence of a diagonal dissection of $P$ into triangles of integer area.
    In fact, utilizing the language introduced in the preceding section, we prove a slightly more general statement that applies to combinatorial polygons.  
    
    Some of the terminology for lattice polygons can be transferred to this context.  
    For example, if $G$ is a combinatorial polygon, then a \emph{diagonal dissection} of $\cp$ is a triangulation of $\cp$ that has no vertices other than the corners.  Also, just as vertices of lattice polygons are colored by parity, the vertices of a combinatorial polygon can be colored with arbitrary colors; a \emph{colored polygon} is a combinatorial polygon with colored vertices.
    If $\cp$ is a colored polygon, then the \emph{boundary word} is the cyclic word consisting of the colors in 
    order\footnote{This is well-defined up to reversing the cyclic ordering; the ambiguity causes no trouble.} around $\cp$.

    Still, to emphasize the transition from lattice polygons to combinatorial ones, we also introduce a bit of new terminology.
    If $\Delta$ is a colored triangle, then we call $\Delta$ \emph{good} if two vertices of $\Delta$ have the same color. 
    Likewise if $\cp$ is a colored polygon, then a \emph{good dissection} of $\cp$ is a diagonal dissection of $\cp$ into good triangles.

    \begin{prop} \label{prop:writemylips}
        A colored polygon $\cp$ has a good dissection if and only if $\word(\cp)$ is contractible.  
    \end{prop}

    A triangle formed by three consecutive vertices of a combinatorial polygon is called {\it external}.
    For example, the dissection in Figure \ref{fig:diagonal} has four external triangles. 
    It is not hard to show that every diagonal dissection includes at least one external triangle.
    

    An external triangle $\Delta$ with (consecutive) vertices $x,y,z$ can be \emph{removed} from a combinatorial $n$-gon $\cp$ by deleting $y$
    and adding an edge connecting $x$ and $z$.  As long as $n\ge 4$ this leaves a combinatorial $(n-1)$-gon 
    $\cp'=\cp\setminus\Delta$.  
    If $G$ is colored, then $\word(\cp')$ is obtained from $\word(\cp)$ by deleting the color of $y$.
    For example, in Figure \ref{fig:diagonal}, removing the $ABA$ triangle results in deleting the letter $B$ from the boundary word.
    Conversely, deleting a letter $X$ from $\word(\cp)$ can be mimicked by removing from $\cp$ an external triangle whose middle vertex corresponds to the deleted $X$.  

    
    \begin{proof}
        Suppose first that $\word(\cp)$ is contractible. This means that there is a sequence of contracting steps that reduces $\word(\cp)$ to a single letter.  Each such step (except the very last one) corresponds to the removal of a good external triangle from the remaining polygon.
        These triangles collectively form a diagonal dissection of $\cp$ into good triangles.  
        
        Conversely, suppose $\cp$ has a diagonal dissection into good triangles.  This dissection includes at least one external triangle.  Removing this triangle executes a contracting step on $\word(\cp)$, and repeating this process results in a word of length two, showing that $\word(\cp)$ is contractible. 
    \end{proof}

    In the context of convex lattice polygons, the Parity Proposition implies that a good dissection corresponds exactly to an integral diagonal dissection.
    
    \begin{cor} \label{cor:readmylips}
        Let $P$ be a convex lattice polygon.  Then $P$ has an integral diagonal dissection if and only if $\word(P)$ is contractible.
    \end{cor}

    \begin{proof}
        Let $\cp$ be a combinatorial polygon with vertices colored so that $\word(\cp) = \word(P)$.  By Proposition \ref{prop:writemylips}, $\word(P)$ is contractible if and only if $\cp$ has a good dissection, which is the case if and only if $P$ has an integral diagonal dissection.
    \end{proof}

\section{You can't go wrong}    \label{sec:diagonal}
    There is a subtlety to the contracting process that we want to address here.  
    Suppose you are given a colored $n$-gon $\cp$, and you want to determine whether or not $\cp$ has a good dissection. 
    If there is no good external triangle, then you know for sure that there is no good dissection. On the other hand, if you find a good external triangle, then you can remove it, with the hope of continuing this process.  However, you may have a choice to make:  if there is more than one good triangle, it seems like you may be at risk of removing the wrong one and heading down a path that does not ultimately lead to a good dissection.
    For example, consider the dodecagon from Figure \ref{fig:polygons}. In Example \ref{ex:diagdiss}, we first removed the good external $ABA$ triangle, leaving an 11-gon\footnote{Despite being Canadian, the first author just learned the word \emph{hendecagon}.}
    $\cp'$, which we then continued to dissect.
    But if we had removed one of the good external $BAB$ triangles instead, we would have been left with a different 11-gon $\cp''$.
    Could it be that $\cp''$ has no good dissection, even though $\cp$ does?

    As it happens, fortunately, you can't go wrong. The following lemma shows that a good external triangle is never a bad option.
    
    \begin{lem}[Diamond lemma, polygon version]\label{lem:diamond}
        Let $\cp$ be a colored $n$-gon, $n\ge 4$. 
        Suppose that $\Delta$ is a good external triangle and let $\cp'=\cp\setminus\Delta$.  
        Then $\cp$ has a good dissection if and only if $\cp'$ does. 
    \end{lem}

    \begin{proof}
        One direction is easy:  if $\cp'$ has a good dissection, then so does $\cp$, since one can add $\Delta$ to a good dissection of $\cp'$ to obtain a good dissection of $\cp$. 
        
        We now assume that $\cp$ has a good dissection and show that $\cp'$ does also. 
        
        We argue by induction on the number of vertices. Consider the base case $n=4$.  If the four vertices of $\cp$ are colored using only two distinct colors, then the triangle $\cp'$ must be good, and we are done. Thus we may assume that at least three distinct colors are used for the vertices of $\cp$.  Since $\cp$ has a good dissection, one can see that there must be two opposite vertices of $\cp$ that are the same color, and $\Delta$ must be a triangle that includes both of them. Hence $\cp'=\cp\setminus\Delta$ also includes these vertices, and therefore $\cp'$ is good.
         
        Now suppose that $n> 4$. Since $\cp$ has a good dissection, there is a good external triangle $\Gamma$ such that $\cp\setminus\Gamma$ has a good dissection. There are some easy cases:  if $\Gamma=\Delta$, then we are done, and if $\Delta\cap\Gamma$ is empty or a single vertex, then since $\cp\setminus\Gamma$ has a good dissection and has fewer vertices than $\cp$, $\cp\setminus (\Gamma\cup\Delta)$ has a good dissection by induction.  Adding back the good triangle $\Delta$, we see that $\cp'$ has a good dissection. 
    
        \begin{figure}[b]
            \centering
            \begin{tikzpicture}[scale=1]
            \input ./diamondlemma.tex
            \end{tikzpicture}
            \caption{Proof of the diamond lemma.}
            \label{fig:diamondlemma}
        \end{figure}
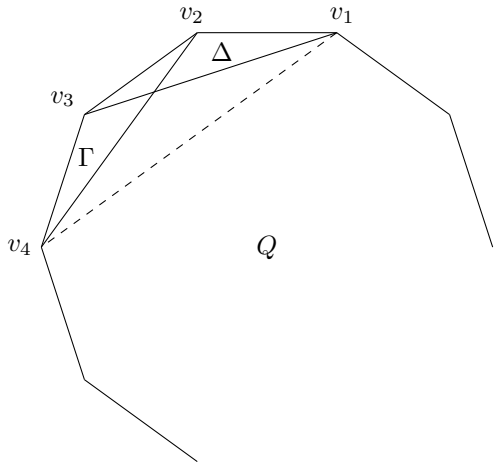
                        
        The remaining possibility is that $\Gamma$ shares an edge with $\Delta$. Label so that $\Delta$ has vertices $v_1,v_2,v_3$ and $\Gamma$ has vertices $v_2,v_3,v_4$, and let $Q$ be the polygon obtained by deleting the quadrilateral $v_1v_2v_3v_4$ from $\cp$.  See Figure \ref{fig:diamondlemma}. Let $v_i$ have color $c_i$, for $i=1,2,3,4$. Suppose first that only two colors are represented here. Then $v_1, v_2, v_4$ is good, so by induction $Q$ has a good dissection, hence adding the triangle $v_1, v_3, v_4$ shows that $\cp'$ has a good dissection.  On the other hand suppose that more than two colors are represented. Then we must have $c_2=c_3$. In this case the boundary word for $\cp\setminus\Gamma$ is the same as the boundary word for $\cp\setminus\Delta=\cp'$. Hence $\cp'$ has a good dissection.
    \end{proof}
    
    The previous lemma is stated in terms of polygons, but it has the following equivalent formulation in terms of words.

    \begin{lem}[Diamond lemma, word version]\label{lem:flakyword}
      Let $w$ be a cyclic word and let $w'$ be obtained from $w$ by a contracting step.  Then $w$ is contractible if and only if $w'$ is contractible.
    \end{lem}

    Readers familiar with combinatorial group theory might recognize elements of the Diamond Lemma for reduction of words in free groups.  Indeed this is exactly what is going on here. In reducing words in free groups, it doesn't matter in what order the reductions are done. Our argument for Lemma \ref{lem:diamond} may be viewed as a geometric proof of this fact. Indeed, many basic facts from combinatorial group theory can be given pictorial proofs of this kind.

\section{A stack of diamonds}\label{sec:alg}
     
    The Diamond Lemma suggests an easy method to decide whether a given cyclic word $w$ is contractible. 
    Given $w$, simply look for any contracting step. As soon as you are able to do a contracting step, do it and continue. 
    If you eventually succeed in contracting $w$ to a single letter, then of course $w$ is contractible; 
    on the other hand if you ever get stuck in the process at a string that has no available contracting steps, then $w$ is not contractible by Lemma \ref{lem:flakyword}. 

    This algorithm has an extremely efficient and simple implementation using a stack. 
    Starting at any letter of $w$, push the letters one by one onto the stack.  After each push, check the top letters of the stack for contracting steps and execute any that are possible. After processing the last letter of $w$, if the stack has a single letter, then you've contracted $w$. Otherwise, Lemma \ref{lem:flakyword} guarantees that $w$ is not contractible. This gives a linear time algorithm for deciding whether the word $w$ is contractible. 
    The integral diagonal dissection shown in Figure \ref{fig:diagonal} is the result of running this algorithm on the dodecagon in Figure \ref{fig:polygons}.
    
    We formalize this observation in the following proposition.
    
    \begin{prop}\label{prop:speedy}
        There is an algorithm that takes as input a cyclic word $w$ of length $n$ and decides in time linear in $n$ whether $w$ is contractible.
    \end{prop}
    
    Note that given a lattice $n$-gon $P$, computing $\word(P)$ also requires only linear time in $n$, as the color of a vertex depends only on the least significant bit of each coordinate.


\section{Poofing} \label{sec:poofing}
    We now shift our attention back to the main question.
    We have shown that if $\word(P)$ is contractible, then $P$ has an integral dissection (a diagonal one, in fact).
    We also know that if $\word(P)$ is not contractible, then $P$ has no \emph{diagonal} integral dissection.
    But we are interested in more than just diagonal dissections.
    If there were no integral dissection at all, how could we prove it?

    A similar problem for rectangles was solved very cleverly by John Thomas \cite{thomas}, using Sperner's Lemma.
    We will give a new version of Sperner's Lemma which allows us to extend Thomas's argument to arbitrary lattice polygons.

\begin{figure}[h]
            \centering
            \begin{tikzpicture}[scale=.5]
            \draw[step=1, gray, very thin] (-.2,-.2) grid (7.2,6.2);

\coordinate (P) at (2,1);
\draw[thick] (P) -- ++(3,0)
    -- ++(1,3)
    -- ++(-3,1)
    -- ++(-2,-2)
    -- cycle;

\draw[thick] (P)++(-1,2) -- +(4,0); 
\draw[thick] (P)++(0,2) -- ++(0,-2) -- ++(1,2) -- ++(0,-2); 
\draw[thick] (P)++(1,2) -- ++(2,-2) -- ++(0,2) -- ++(-2,2); 
\draw[thick] (P)++(3,2) -- ++(1,1);

\draw[fill=black] (P) circle (2pt);
\draw[fill=black] (P)+(3,0) circle (2pt);
\draw[fill=black] (P)+(4,3) circle (2pt);
\draw[fill=black] (P)+(1,4) circle (2pt);
\draw[fill=black] (P)+(-1,2) circle (2pt);

\draw[fill=black] (P)+(0,2) circle (2pt);
\draw[fill=black] (P)+(1,2) circle (2pt);
\draw[fill=black] (P)+(3,2) circle (2pt);
\draw[fill=black] (P)+(1,0) circle (2pt);

\coordinate (1) at (14,2);
\coordinate (2) at (16,2);
\coordinate (3) at (17,4.5);
\coordinate (4) at (15,6);
\coordinate (5) at (13,4.5);
\coordinate (6) at (13.8,4.1);
\coordinate (7) at (15.2,3.8);
\coordinate (8) at (16,4.2);
\coordinate (9) at (15,2.2);

\draw (1)--(2)--(3)--(4)--(5)--cycle;
\draw (1)--(6);
\draw (5)--(6);
\draw (6)--(7);
\draw (1)--(7);
\draw (9)--(7);
\draw (1)--(9);
\draw (9)--(2);
\draw (7)--(2);
\draw (7)--(8);
\draw (2)--(8);
\draw (8)--(3);
\draw (8)--(4);
\draw (5)--(8);

\filldraw[fill=white!80!black] (1)--(2)--(9)--cycle;
\filldraw[fill=white!80!black] (5)--(6)--(7)--(8)--cycle;

\draw[dashed] (6)--(8);

\draw[->] (12,5) .. controls (10,6) and (9,4) .. (8,4);
            \end{tikzpicture}
            \caption{Poofing up a dissection into a triangulation.}
            \label{fig:poofing}
        \end{figure}
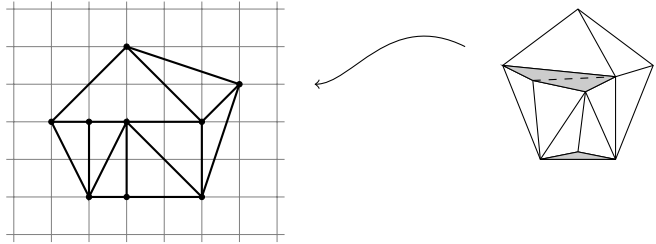

    Typically, Sperner's Lemma is proved and applied in the context of triangulations of combinatorial $n$-gons (often $n=3$).
    A concern one might have about Definition \ref{def:triangulation} is that at first glance, it seems that these objects are not as general as dissections.
    For example, Figure \ref{fig:dissection} does not appear to depict a triangulation of a pentagon.  This is because in any triangulation, two triangles (viewed officially as triples of vertices) intersect in either the empty set, a single vertex, or two vertices.  Geometrically, this means that two triangles of $|T|$ that intersect in more than a vertex must meet along an entire edge of both triangles.  This is not the case in Figure \ref{fig:dissection} (reproduced in Figure \ref{fig:poofing}), where the top left triangle intersects the triangles below in only part of an edge.  So, maybe Sperner's Lemma will not apply to dissections.

    Luckily, there is a way around this.          
    The idea is to view a dissection $D$ as a \emph{picture} of a triangulation $T$, specifically a picture in which some triangles may be drawn with area zero.  For example, the bottom edge of the dissected pentagon in Figure \ref{fig:poofing} contains a third vertex of the dissection.  We can view all three vertices as forming a triangle that happens to have area zero in this particular planar picture.

        Generally where more than two vertices of the dissection $D$ are collinear, we can view the collinear points as degenerate polygons which ``poof up'' the dissection.  Using the triangles of $D$ together with the ``poofagons'' that are invisible in $D$, we can build a triangulation $T$ of which $D$ is a partially degenerate picture. This is illustrated in Figure \ref{fig:poofing}, in which the poofagons are shaded.

        The horizontal segment running through the middle of the pentagon contains not just three but four vertices of the dissection.  The resulting poofagon is a quadrilateral, so we must add an extra edge to triangulate it.  This edge is shown as a dashed line in the triangulation.  We made a choice here:  we could have chosen the other diagonal of this quadrilateral instead.  So the triangulation $T$ is not completely well-defined, but this doesn't matter since in the dissection we care about, all those triangles are degenerate.

        Note also that no new vertices are formed when creating $T$ from $D$.
        An eager reader can work out a general procedure for carrying out this process; alternatively, details can be found in \cite{chapter3}.  For our purposes it is enough to record the following fact.

        \begin{fact}
            If $P$ is a dissected polygon in the plane, then there is a triangulation $T$ of a combinatorial polygon and a map $f:|T|\to P$ so that $f$ maps each triangle of $|T|$ either to a triangle of the dissection or to a degenerate triangle in $\R^2$.  Moreover $f$ gives a one-to-one correspondence from the vertices of $T$ to the vertices of the dissection.
        \end{fact}

        Sperner's Lemma is back in play!

\section{Non-abelian Sperner} \label{sec:nonabelian}

    Sperner's original lemma \cite{sperner,su} is about triangulations of a triangle, with vertices colored using three colors. Many versions of Sperner's Lemma, including the original, identify a property of a coloring on the boundary that guarantees the existence of a triangle somewhere in the triangulation whose vertices are all different colors.  Our version does this too.

    \begin{defn}[Colored triangulation]
        Let $\cp$ be a combinatorial polygon.  A \emph{colored triangulation} of $\cp$ is a triangulation $T$ of $\cp$ and an assignment of a color to each vertex of $T$.  If $\cp$ is colored to begin with, then the coloring of $T$ is required to extend the coloring of $\cp$.  The \emph{boundary word} $\word(T)=\word(\cp)$ is the cyclic word consisting of the colors of the corners (i.e., of $\cp$) in cyclic order.
    \end{defn}
    
    \begin{defn}[Tricolor]
        A colored triangle is called \emph{tricolor} if its vertices all have different colors.
    \end{defn}

    \begin{thm}[Non-abelian Sperner]\label{thm:nonabelian}
        Let $T$ be a colored triangulation of a combinatorial polygon. 
        If $\word(T)$ is not contractible, then $T$ contains a tricolor triangle. 
    \end{thm}

    \begin{proof}
        Assume $T$ is a triangulation of a combinatorial $n$-gon whose boundary word is not contractible.  We will show that there is a tricolor triangle. We proceed by induction on the number of triangles in the triangulation. If there is just one, then $n=3$. Since $\word(T)$ is not contractible, the corners of $T$ must all be different colors, so $T$ is already the tricolor triangle we seek.
        
        Now suppose $T$ has some number of triangles. Choose a boundary edge $uv$ and let $x$ be the other vertex of the triangle containing $uv$.  If $uvx$ is tricolor, we are done, so suppose this is not the case.  Our strategy will be to delete the triangle $uvx$ and proceed by induction.  There are three possibilities for the location of $x$: (1) in the interior of $T$; (2) a corner of $T$, adjacent to $u$ or $v$; or (3) a corner of $T$, not adjacent to $u$ or $v$.
        See Figure \ref{fig:spernerproof}.
        
        \begin{figure}[h]
            \centering
            \begin{tikzpicture}[scale=1]
            \coordinate (v) at (-5,0);
\coordinate (u) at (-5,.8);
\coordinate (w) at (-4.4,.5);
\coordinate (x) at (-4.4,1.3);
\coordinate (y) at (-4.4,-.3);

\draw (y)--(v)--(u)--(x);
\draw (v)--(w)--(u);

\draw (u) node[anchor=east]{$u$};
\draw (v) node[anchor=east]{$v$};
\draw (w) node[anchor=west]{$x$};

\draw (-4.7,-1) node[below]  {Case 1};

\begin{scope}[xshift=-1cm]
\coordinate (u1) at (0,0);
\coordinate (v1) at (0,.8);
\coordinate (w1) at (.6,1.3);
\coordinate (x1) at (1,1.3);
\coordinate (y1) at (.6,-.3);

\draw (y1)--(u1)--(v1)--(w1)--(x1);
\draw (u1)--(w1);

\draw (u1) node[anchor=east]{$v$};
\draw (v1) node[anchor=east]{$u$};
\draw (w1) node[anchor=south ]{$x$};

\draw (.5,-1) node[below]  {Case 2};
\end{scope}

\begin{scope}[xshift=8cm]
\coordinate (v) at (-5,0);
\coordinate (u) at (-5,.8);
\coordinate (w) at (-3,.5);
\coordinate (x) at (-4.6,1.3);
\coordinate (y) at (-4.6,-.3);
\draw (w) -- +(-.1,.5);
\draw (w) -- +(-.1,-.5);
\draw (y)--(v)--(u)--(x);
\draw (u)--(w)--(v);

\draw (u) node[anchor=east]{$u$};
\draw (v) node[anchor=east]{$v$};
\draw (w) node[anchor=west]{$x$};

\draw (-4,-1) node[below]  {Case 3};
\end{scope}
            \end{tikzpicture}
            \caption{Where's $x$?}
            \label{fig:spernerproof}
        \end{figure}
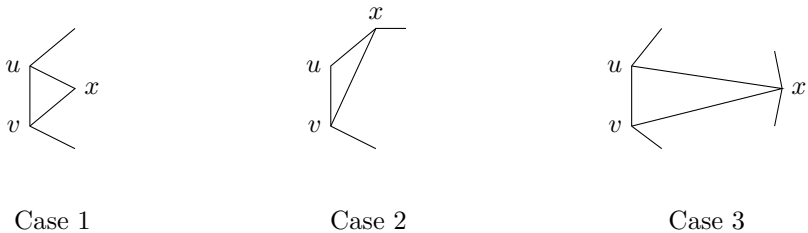

        It is an exercise to show that it is always possible to choose the edge $uv$ so that we are in Case 1 or Case 2.  (Alternatively one could deal with Case 3 directly, without much trouble.)  In each of these cases, deleting the triangle $uvx$ from $T$ results in a triangulation $T'$ with fewer triangles than $T$. And in each case, $\word(T)$ and $\word(T')$ differ by a single contracting step. Lemma \ref{lem:flakyword} implies in each case that $\word(T')$ is not contractible.  Therefore, by induction, $T'$ contains a tricolor triangle, which is also a tricolor triangle of $T$.  This completes the proof.
    \end{proof}

    If we view the Non-abelian Sperner Lemma as an assertion about colored polygons (rather than triangulations), then we can state it as a biconditional.  Compare with Proposition \ref{prop:writemylips}.

    \begin{thm}[Non-abelian Sperner]\label{thm:nonabelian2}
        Let $\cp$ be a colored polygon. Then there exists a colored triangulation of $\cp$ with no tricolor triangles
        if and only if $\word(\cp)$ is contractible. 
        
        Moreover, if such a triangulation exists, then one exists with no interior vertices.
        %
        %
    \end{thm}

\begin{figure}[h]
        \centering
        \begin{tikzpicture}[scale=.5,
            vertex/.style={circle,
                              inner sep=0pt,minimum size=2pt}]
        \foreach \x in {0,...,10}
    \coordinate (\x) at (36*\x : 3);

\foreach \x/\y in {0/1,1/2,2/3,3/4,4/5,5/6,6/7,7/8,8/9,9/10}
\draw (\x) -- (\y);
\draw (0) node [anchor=180]{$A$};
\draw (1) node [anchor=216]{$B$};
\draw (2) node [anchor=252]{$C$};
\draw (3) node [anchor=288]{$D$};
\draw (4) node [anchor=324]{$A$};
\draw (5) node [anchor=0]{$C$};
\draw (6) node [anchor=36]{$B$};
\draw (7) node [anchor=72]{$A$};
\draw (8) node [anchor=108]{$D$};
\draw (9) node [anchor=144]{$C$};
        \end{tikzpicture}
        \caption{This colored decagon cannot be triangulated and colored without a tricolor triangle.}
        \label{fig:decagon}
    \end{figure}
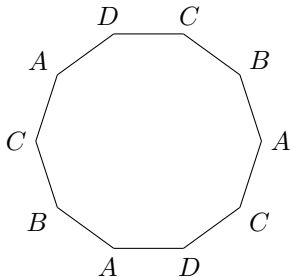

    \begin{example}\label{ex:nonabsper}
        Figure \ref{fig:decagon} shows a combinatorial decagon whose boundary word is $(ABCDACBADC)$, which is also the boundary word of the lattice decagon from Figure \ref{fig:polygons}. This word is not contractible, as there are no contracting steps available.  By the Non-abelian Sperner Lemma, any combinatorial triangulation of the decagon with any coloring of the interior vertices will contain a tricolor triangle.
    \end{example}
    
\section{Impossibility} \label{sec:impossibility}

    Here is the heart of our main theorem.
    
    \begin{thm}\label{thm:integerclassic} 
        Let $P$ be a lattice polygon, with corners colored according to their parities.  
        If $\word(P)$ is not contractible, then $P$ has no integral dissection.
    \end{thm}
    
    \begin{proof}
        Assume $\word(P)$ is not contractible.  Suppose that $P$ is dissected into lattice triangles, and as usual color each vertex of the dissection according to its parity.
        
        Now form a poofed triangulation $T$ with corresponding map $f:|T|\to P$ as described in Section \ref{sec:poofing}.  
        In $T$, assign each vertex the same color as the corresponding vertex of the dissection.  
        
        Note that $T$ satisfies the hypothesis of Theorem \ref{thm:nonabelian}.  
        Therefore $T$ contains a tricolor triangle $\Delta$.  What happens to $\Delta$ when we apply the map $f$?  By the Collinear Corollary, $f(\Delta)$ cannot be a degenerate triangle, so it must be one of the triangles of the dissection.  Since it is tricolor, its area is not an integer.
    \end{proof}
    
    \begin{example} \label{ex:impossibility}
        The $3\times 5$ rectangle (Figure \ref{fig:polygons}) has boundary word $(ABCD)$, and therefore cannot be dissected into integer area lattice triangles.
        This is an instance of the theorem proved by Thomas \cite{thomas}, and it can be shown using the usual Sperner Lemma.
        
        The hexagon in Figure \ref{fig:polygons} has boundary word $(ABCABC)$, and therefore also has no integral dissection.
        In this case we could deduce this using another classical version of Sperner's Lemma in which orientations are considered. 



        To the best of our knowledge, though, no previously known version of Sperner's Lemma can treat the decagon from Figure \ref{fig:polygons}.
        Its boundary word $(ABCDACBADC)$ is not contractible (see Example \ref{ex:nonabsper}), so by Theorem \ref{thm:integerclassic}, it has no integral dissection.

        And lest we forget, the dodecagon from Figure \ref{fig:polygons} has contractible boundary word, and we have explicitly seen that it has an integral (diagonal) dissection.
    \end{example}
    
    Our route to Theorem \ref{thm:integerclassic} started with a careful study of Thomas's ideas. 
    He essentially used Sperner's Lemma to provide a significant first step in solving a dissection problem that he and Fred Richman posed \cite{richman}.  
    Shortly thereafter the problem was fully resolved by Paul Monsky \cite{monsky,thebook,steinbook}, who proved that a square cannot be dissected into an odd number of triangles of equal area.  Monsky's theorem and its clever proof have captured the imagination of many mathematicians and inspired quite a bit of new and fascinating mathematics.
    Attempting to better understand Monsky's work led us to Thomas, and then to our main theorem.

\section{Altering integral triangle} \label{sec:appendix}

    We finally address the triangle case of our main question, which establishes the equivalence of (1) and (2) in Theorem \ref{thm:stacked}.
        
    \begin{thm}[The triangle case] \label{thm:trianglecase}
        Every lattice triangle $P$ with (positive) integer area can be dissected into lattice triangles of area 1.
    \end{thm}

    To prove this, we first put our triangle in a standard position. This idea is used in many proofs of Pick's Theorem and other contexts. See, for example, the $x$-axis Lemma in \cite{rabinowitz}.
    Here we will make use of linear transformations defined by the $2\times 2$ matrices in the group 
    \[
    GL_2(\Z)= \{ A\in M_{2,2}(\Z) \ \ | \ \ \det A = \pm 1 \}.
    \]
    If $A\in GL_2(\Z)$ then $A$ has an inverse whose entries are also integers, and both $A$ and $A^{-1}$ preserve the lattice. Such automorphisms of the lattice preserve straight lines and areas, and hence also preserve dissections into integer area triangles.

    \begin{lem}
        Let $P$ be a lattice triangle.  Then using  matrices in $GL_2(\Z)$ and lattice translations, we can map $P$ to a triangle $P'$ with vertices at $(0,0), (d,0), (p,q)$, where $d, p ,q$ are integers such that $d>0$,  and $1\le p \le q$.
    \end{lem}

    \begin{proof}
        First a lattice translation can be used to move one vertex of the triangle to the origin. Suppose the other two vertices map to the lattice points $v_1$ and $v_2$, labeled so that $\det(v_1,v_2)>0$. Let $(a,b)$ be the coordinates of $v_1$  and let $d=\gcd(a,b)$. Write $d = ra + sb$ as an integer linear combination of $a$ and $b$. Then the matrix
        \[
        \begin{pmatrix}
                r & s \\
                -b/d & a/d 
            \end{pmatrix}
        \]
        has determinant $1$ and maps $v_1$ to $(d,0)$ and $v_2$ maps to a  point $(t,q)$ with $q>0$. 
        
        
        Finally, the integer $t$ can be reduced mod $q$. Specifically, there exists an integer $p$ with $1 \leq p \leq q$ such that $p=t+kq$, where $k\in \Z$. The matrix
        $$ \begin{pmatrix}
                1 & k \\
                0 & 1 
                \end{pmatrix}
        $$
        is in $GL_2(\Z)$, maps $(t,q)$ to $(t+kq, q)=(p,q)$, and fixes the other two vertices $(0,0)$ and $(d,0)$. Thus we are done. 
    \end{proof}

    Note that any such map preserves equality of colors.  (Recall that all lattice points are colored by parity.)
    That is, if two points have the same color, then their images also have the same color. To see this, note that two points have the same color if and only if the midpoint of the segment joining them is a lattice point. 
    
    Observe also that if $P$ is any triangle and $x$ is any point of the interior of $P$, then we may dissect $P$ into three triangles by connecting $x$ to the three corners of $P$. Similarly, if $x$ is in the interior of an edge of $P$, we may dissect $P$ into two triangles, each having $x$ as a vertex.  In either case, we refer to this as ``using $x$ to dissect $P$."  Also note that if $P$ has two vertices of the same color and $x$ has the same color as any vertex of $P$, then no triangle of this dissection will have three vertices of different colors. Hence, all triangles of this dissection will have integer area.

    \begin{proof}[Proof of Theorem \ref{thm:trianglecase}]
        We will proceed by induction.  If the area of $P$ is 1, then nothing needs to be done. So suppose that the area of $P$ is larger than $1$.  By induction, it suffices to dissect $P$ into lattice triangles of smaller integer area.        
        
        Since the area of $P$ is an integer, there are two vertices that are the same color. Label these $v_0, v_1$, and the third vertex $v_2$, so that the signed area of $v_0v_1v_2$ is positive.  By the previous lemma, there is a matrix in $GL_2(\Z)$ that takes $v_0, v_1, v_2$ to a triangle $P'$ with vertices $(0,0), (d,0), (p,q)$, respectively. Since $(0,0)$ and $(d,0)$ are the same color, $d$ must be even.  If $d>2$, then we may use $(2,0)$ to dissect $P'$ into two triangles of smaller integer area. Thus we may assume $d=2$. If $q$ is even, then we may use $(1,0)$ to dissect $P'$ into two triangles of area $q/2$. Also note that since $P$ and $P'$ have the same area, we may assume $q>1$.
         
        The remaining case is that $P'$ has vertices $(0,0)$, $(2,0)$ and $(p,q)$, where $q>1$ is odd and $1\leq p \le q$. If $p$ is also odd, then $(p,q)$ has the same color as $(1,1)$, which lies in $P'$ and may be used to dissect $P'$. If $p$ is even, then $(p,q)$ has the same color as $(2,1)$, which lies in $P'$ and may be used to dissect $P'$.
         
        In all cases, $P'$ and therefore $P$ can be dissected into lattice triangles of smaller integer area. 
    \end{proof}

\section{Epilogue}
    
    \subsection*{What's new?} 
        Dissection problems are not new, nor is Sperner's Lemma, nor is the connection between the two, as we have discussed.  The concept of contractibility is not new.  The core idea underlying the proof of our main theorem, namely that \emph{colorings define maps} (see below), is also not new.  For instance there is a very nice proof of Sperner's Lemma based on this idea due to the economists A.~McLennan and R.~Tourkey \cite{econsperner}, which we learned from R.~Schwartz \cite{schwartzsperner}. Even the specific idea of using the parity coloring of lattice points to define a map to $K_4$ was used by D.~Rudenko \cite{rudenko}, whose Lemma 3 implies a weakened version of our Theorem \ref{thm:integerclassic}. (The additional power of our results comes from using homotopy instead of homology.)

        It seems a little strange that our theorems have gone undiscovered for the last fifty years.
        One possible explanation is that in the study of polytopes, ``integral'' and ``lattice'' are typically synonymous; thus the volume of an integral (in that sense) $d$-dimensional polytope is an integer multiple of $1/d!$.  To our knowledge, little has been written about simplices whose volumes are actual integers.

    \subsection*{Colorings define maps} 
        
        We believe the most natural way to view Theorems \ref{thm:nonabelian} and \ref{thm:nonabelian2} is through a topological lens.
        The key idea is that colorings define maps.
        A colored polygon $\cp$ using $m$ colors naturally defines a map from $\cp$ to a complete graph $K_m$ with nodes labeled by the colors:  each vertex of $\cp$ is mapped to the vertex of $K_m$ of the corresponding color, and each edge is mapped in the obvious way either to an edge of $K_m$ or to a point if both endpoints have the same color.  Topologically, $\cp$ is a circle so this map is a loop $f$ in $K_m$.  The combinatorial notion of contractibility for words (Definition \ref{def:contract}) is meant to mimic the topological contractibility of loops in a graph.  In particular, the word $\word(\cp)$ is contractible if and only if the loop $f$ is (topologically) contractible in the space $K_m$.
        
        If $T$ is a colored triangulation of $\cp$ (using no new colors), then this extended coloring defines an extended map $F$ from the disk $|T|$ to the 2-dimensional simplicial complex $X$ obtained from $K_m$ by adding all possible triangles.  (This complex is the 2-skeleton of an $(m-1)$-dimensional simplex.)  The map is defined as above on the vertices and edges, and each triangle is then mapped by extending the already-defined map on its boundary.  The image of a triangle is a vertex, edge, or triangle, depending on the number of distinct colors of its vertices.

        In particular if $T$ is a colored triangulation of a polygon using $m$ colors then the loop $f:G\to K_m$ is contractible in $X$, since it bounds the disk defined by $F$.  The loop $f$ may or may not be contractible in $K_m$ itself; the point is that $f$ is indeed contractible in $K_m$ if the entire image of $F$ is contained in $K_m$, which happens if and only if $T$ has no tricolor triangles.

        
        \begin{proof}[Topological proof of Non-abelian Sperner] 
            The given coloring defines a map $F$ from $|T|$ to $X$, the 2-skeleton of an $(m-1)$-simplex (where $m$ is the number of colors in the coloring). If $T$ contains no tricolor triangle, then by the observation above the image of $F$ is contained in the $1$-skeleton of $X$, which is the graph $K_m$.  Thus the boundary loop is contractible in $K_m$ (not just in $X$), and therefore the boundary word is a contractible word.
        \end{proof}


    \subsection*{Why ``non-abelian"?}
        Our version of Sperner's Lemma arises from the interpretation of boundary words as loops in the graph $K_4$, whose fundamental group is a non-abelian free group. The hypothesis is that the boundary word represents a non-trivial element of this group.
        
        If we had only three colors, as in the usual Sperner Lemma, the boundary word would represent a loop in the space $K_3$, whose fundamental group is (isomorphic to) the abelian group $\Z$.  The hypotheses in the usual versions of Sperner's Lemma also express that the boundary is non-trivial, but on a circle, the non-contractibility of a loop, or equivalently the property of having non-zero winding number, is detected more simply, for instance by counting the (signed) number of $AB$ edges. 
        
        In the fundamental group of $K_4$ this doesn't work.  The decagon in Figure \ref{fig:polygons} has boundary word $(ABCDACBADC)$, which has a signed total of zero $AB$ edges (and any other type of edge), and yet is not contractible.  This word represents a commutator in the fundamental group, so its non-triviality is not detected by any abelian measurement.

    \subsection*{Convexity}
        Almost everything about our theorems is combinatorial.  
        In particular, the convexity hypothesis is not used in the application of the Non-abelian Sperner's Lemma.
        Any lattice polygon, convex or not, with a non-contractible boundary word has no integral dissection.
        However, convexity \emph{is} used in the construction of diagonal dissections:
        we use the fact that we can remove any external triangle we want.

        \begin{figure}[h]
            \centering
            \begin{tikzpicture}[scale=.4,rotate=90]

\draw[step=1, gray, very thin] (-.2,-.2) grid (4.2,4.2);

\coordinate (P) at (2,2);
\draw[very thick] (P) -- ++(-1,-1);
\draw[very thick] (P)++(-1,-1) -- ++(2,1) -- ++(-2,1);
\draw[very thick] (P)-- ++(-1,1) ;
\draw[fill=white] (P) circle (2pt);
\draw[fill=white] (P)+(-1,-1) circle (2pt);
\draw[fill=white] (P)+(1,0) circle (2pt);
\draw[fill=black] (P)+(-1,1) circle (2pt);

\draw[thick,dashed] (P)++(-1,1) -- ++(0,-2);
            \end{tikzpicture}
            \caption{The boundary word of this polygon, $(ACAB)$, is contractible.}
            \label{fig:nonconvex}
        \end{figure}
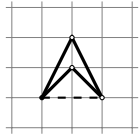
        
        It turns out that this poses no problem in the non-convex case either, if we are willing to use triangles with \emph{negative} (but still integer) areas!  (See Equation \eqref{eq:area}.)  For example, consider the dart-shaped quadrilateral in Figure \ref{fig:nonconvex}. This polygon has no dissection into lattice triangles of integer area, in the usual sense of the word ``dissection.'' Yet its boundary word $(ACAB)$ is contractible.
        The recipe in the proof says to add the ``diagonal'' joining the two $A$ vertices, which is a horizontal segment shown in the figure as a dashed line. 
        This segment lies outside the polygon.
        Nevertheless, if we add it anyway, it ``cuts'' the polygon into two triangles, one of area $2$ and one of area $-1$. We might think of this as a sort of ``signed'' dissection.
        
        We do not know a general criterion for a non-convex lattice polygon to have a genuine integral dissection.

    \subsection*{Higher dimensions}
        Sperner's Lemma and many of its generalizations can be applied in all dimensions.
        Our version, as stated, is about contractible loops and the disks they bound.  It differs from Sperner's original lemma by allowing an arbitrary number of colors, even in just two dimensions.  This is what necessitates the contractibility hypothesis.
        It came as a pleasant surprise to us that contractibility in $K_4$, interpreted in the context of lattice polygons and the parity coloring, corresponds exactly to tiling lattice polygons with integer area triangles.

        In higher dimensions, the colorings-define-maps principle readily applies.  We describe one related theorem; see also \cite{musin}.
        
        Let $P$ be a simplicial $n$-dimensional polytope (this means that each $(n-1)$-dimensional face of $P$ is combinatorially a simplex), with $n\ge 3$.  Suppose $P$ is triangulated and its vertices are colored. 
        If the signed count of boundary simplices with colors $C_1,\ldots,C_n$ is nonzero then there must be a simplex of the triangulation that uses $n+1$ distinct colors.  (Defining the signed count is slightly tricky because it requires the notion of \emph{orientation}.  One could instead ignore the signs and require an odd number of boundary simplices colored $C_1,\ldots,C_n$.)  
        
        Notice, though, that other than having an arbitrary number of colors, this looks a lot like the usual $n$-dimensional Sperner Lemma.  Why?  Here the coloring defines a map from the boundary of $P$ into a space $X$ which is the $(n-1)$-skeleton of a big simplex.  This map represents an element of the homotopy group $\pi_{n-1}(X)$.  By an elementary algebraic topology argument, this map is contractible if and only if all signed counts of oriented simplices come out to zero; there is nothing non-abelian here. 
        
        Unfortunately we don't know how to use these observations to say anything interesting about volumes or dissections of higher dimensional lattice polytopes.

    \subsection*{Surface topology}
        Consider a combinatorial hexagon $H$ with non-contractible boundary word $w=(ABCABC)$.  This object was realized as a lattice polygon in Figure \ref{fig:polygons} (upper right).  The loop in $K_3$ associated to this word wraps around the circle twice.  There is no continuous map from a disk to a circle that wraps the boundary of the disk twice around the circle, so any triangulated disk with boundary word $w$ and with colored vertices must contain a tricolor triangle.  This is Sperner's Lemma.
    
\begin{figure}[h]
            \centering
            \begin{tikzpicture}[scale=.4]
            \coordinate (P) at (0,0);
\coordinate (Q) at (2,0);
\coordinate (R) at (4,0);
\coordinate (S) at (.5,3);
\coordinate (T) at (2.5,3);
\coordinate (U) at (4.5,3);

\draw (P) node[anchor=north]{$A$};
\draw (Q) node[anchor=north]{$B$};
\draw (R) node[anchor=north]{$C$};
\draw (S) node[anchor=south]{$A$};
\draw (T) node[anchor=south]{$B$};
\draw (U) node[anchor=south]{$C$};

\draw (-1,0) -- (6,0);
\draw (6,3)   -- (-1,3);
    decorate {-- (-1,0)};

\draw (P) -- (S) -- (Q) -- (T) -- (R) -- (U);

\draw (U) -- (7,0);
\draw[white,fill=white,decoration={zigzag,segment length=8pt}]  
    (6,0)
    decorate {--  (6,3)}
    -- (7.1,3) -- (7.1,-.1) -- (6,0);
\draw (S) -- (-3,0);
\draw[white,fill=white,decoration={zigzag,segment length=8pt}]  
    (-1,3)
    decorate {--  (-1,0)}
    -- (-3.1,-.1) -- (-3.1,3) -- (-1,3);

\draw[decorate,decoration={zigzag,segment length=8pt}] 
    (6,0) -- (6,3);
\draw[decorate,decoration={zigzag,segment length=8pt}] 
    (-1,0) -- (-1,3);

\draw[->] (8,0) -- (8,2);
\draw     (8,2) -- (8,3);
\draw[->] (-3,3) -- (-3,1);
\draw     (-3,1) -- (-3,0);
            \end{tikzpicture}
            \caption{The boundary word $(ABCABC)$ can be triangulated with no tricolor triangles on a M\"obius strip.}
            \label{fig:mobius}
        \end{figure}
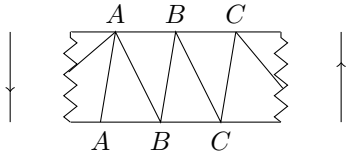

        On the other hand, there \emph{does} exist a different surface $M$ with one boundary component with the property that there is a map from $M$ to a circle that wraps the boundary of $M$ around the circle twice.  This surface is the M\"obius strip.  And indeed, if we view $H$ as the boundary of a M\"obius strip rather than the boundary of a disk, then it is possible to triangulate $H$ with no tricolor triangle.  See Figure \ref{fig:mobius}.
    
        Another way to say this is that a small hexagon on the projective plane cuts the surface into a disk and a M\"obius strip.  If the hexagon is colored with word $(ABCABC)$ then it can be triangulated in the projective plane without tricolor triangles, but only using the M\"obius strip and not the disk.
            
        Similarly, our familiar decagon labeled by the commutator $(ABCDACBADC)$ (also from Figure \ref{fig:polygons}) can be triangulated without tricolor triangles if it is placed on a small region of a torus.  Again one has to use the ``outside,'' i.e., the punctured torus, rather than the disk, to construct such a triangulation.  It is an oddly satisfying exercise to draw such a triangulation.
        
        One can construct higher genus examples as well, starting from presentations of the fundamental group of the corresponding surface.

\subsection*{Acknowledgements}
    Abrams thanks the Distinguished Guest Scientist Fellowship Programme of the Magyar Tudom\'anyos Akad\'emia for their support in 2022.
    Pommersheim thanks the Fulbright U.S.~Scholar Program for their support. Both authors are deeply grateful to Dezs\H{o} Mikl\'os and Andr\'as Stipsicz for their roles in making this work possible.

\bibliographystyle{alphaurl}
\bibliography{Refs.bib}

\end{document}